# CHOQUET EXPECTATION AND PENG'S *G*-EXPECTATION

By Zengjing Chen[1], Tao Chen and Matt Davison[2]

*Shandong University and University of Western Ontario*

In this paper we consider two ways to generalize the mathematical expectation of a random variable, the Choquet expectation and Peng's *g*-expectation. An open question has been, after making suitable restrictions to the class of random variables acted on by the Choquet expectation, for what class of expectation do these two definitions coincide? In this paper we provide a necessary and sufficient condition which proves that the only expectation which lies in both classes is the traditional linear expectation. This settles another open question about whether Choquet expectation may be used to obtain Monte Carlo-like solution of nonlinear PDE: It cannot, except for some very special cases.

**1. Introduction.** The concept of expectation is clearly very important in probability theory. Expectation is usually defined via

$$E\xi = \int_{-\infty}^{\infty} x \, dF(x),$$

where $F(x) := P(\xi \leq x)$ is the distribution of random variable $\xi$ with respect to the probability measure $P$. Alternatively, the expectation $E\xi$ can be written as

$$(0.1) \qquad E\xi = \int_{-\infty}^{0} [P(\xi \geq t) - 1] \, dt + \int_{0}^{+\infty} P(\xi \geq t) \, dt,$$

which implies the relation between mathematical expectation and probability measure. One of the properties of mathematical expectation is its

Received November 2003; revised June 2004.

[1]Supported by the Fields Institute, NSFC (10325106), FANEDD (200118) and NSF of China 10131030.

[2]Supported in part by grants from the Natural Sciences and Engineering Research Council of Canada and MITACS.

linearity: for given random variables $\xi$ and $\eta$,

$$(0.2) \qquad E(\xi + \eta) = E\xi + E\eta.$$

This is equivalent to the additivity of probability measure, that is,

$$(0.3) \qquad P(A + B) = P(A) + P(B) \qquad \text{if } A \cap B = \varnothing.$$

From this viewpoint, we sometimes call mathematical expectation (resp. probability measure) linear mathematical expectation (resp. linear probability measure). It is easy to define conditional expectation using the additivity of mathematical expectations, that is, the conditional expectation $\eta$ of a random variable $\xi$ under a given $\sigma$-field $\mathcal{F}$ is a $\mathcal{F}$-measurable random variable such that

$$(0.4) \qquad E\xi I_A = E\eta I_A \qquad \forall A \in \mathcal{F}.$$

It is well known that linear mathematical expectation is a powerful tool for dealing with stochastic phenomena. However, there are many uncertain phenomena which are not easily modeled using linear mathematical expectations. Economists have found that linear mathematical expectations result in the Allais paradox and the Ellsberg paradox, see Allais (1953) and Ellsberg (1961). How to deal with uncertain phenomena which cannot be well explained by linear mathematical expectations? One idea is to examine nonlinear expectations. How should nonlinear mathematical expectations be defined? Choquet (1953) extended the probability measure $P$ in (0.1) to a nonlinear probability measure $V$ (also called the capacity) and obtained the following definition $C(\xi)$ of nonlinear mathematical expectations (called the Choquet expectation):

$$C(\xi) := \int_{-\infty}^{0} [V(\xi \geq t) - 1] \, dt + \int_{0}^{+\infty} V(\xi \geq t) \, dt.$$

Because $V$ no longer has property (0.3), the above Choquet expectation $C(\xi)$ usually no longer has property (0.2). Choquet expectations have many applications in statistics, economics, finance and physics. Unfortunately, scientists also find that it is difficult to define conditional Choquet expectations in terms of Choquet expectations. Many papers study the Choquet expectation and its applications, see, for example, Anger (1977), Dellacherie (1970), Graf (1980), Sarin and Wakker (1992), Schmeidler (1989), Wakker (2001), Wasserman and Kadane (1990) and the references therein. Peng (1997, 1999) introduced a kind of nonlinear expectation (he calls it the $g$-expectation) via a particular nonlinear backward stochastic differential equation (BSDE for short). Using Peng's $g$-expectation, it is easy to define conditional expectations in the same way as in (0.4). Some applications of Peng's $g$-expectation in economics are considered in Chen and Epstein (2002). An open question raised by Peng is the following: What is the relation between Choquet expectation and Peng's $g$-expectation? We note that Peng's $g$-expectations can be defined only in a BSDE framework, while Choquet capacities and expectations makes sense in more general settings. In this paper when we compare the two objects, we do so after making suitable restrictions to Choquet expectations. That said, does there exist a Choquet expectation whose restriction to the domain $L^2(\Omega, \mathcal{F}, P)$ of $g$-expectation is equal to a $g$-expectation?

An earlier work by Chen, Kulperger and Sulem (2002) shows that the answer is yes for certain special random variables. In this paper we shall further study this question and obtain a necessary and sufficient condition for this open question. This settles another open question about whether Choquet expectation may be used to obtain Monte Carlo-like solution of nonlinear partial differential equations (PDE): It cannot, except for some very special cases.

**2. Notation and lemmas.** In this section we introduce the concepts of Choquet expectation and Peng's $g$-expectation. For convenience, we include some related lemmas used in this paper.

*Capacity and Choquet expectation.* We now introduce the concepts of capacity and Choquet expectation.

DEFINITION 1.    1. Random variables $\xi$ and $\eta$ are called comonotonic if

$$[\xi(\omega) - \xi(\omega')][\eta(\omega) - \eta(\omega')] \geq 0 \qquad \forall \omega, \omega' \in \Omega.$$

2. (*Comonotonic additivity*). A real functional $F$ on $L^2(\Omega, \mathcal{F}, P)$ is called comonotonic additive if

$$F(\xi + \eta) = F(\xi) + F(\eta) \qquad \text{whenever } \xi \text{ and } \eta \text{ are comonotonic.}$$

3. A set function $V : \mathcal{F} \longrightarrow [0, 1]$ is called a capacity if:
   (i) $V(\varnothing) = 0, V(\Omega) = 1$;
   (ii) If $A, B \in \mathcal{F}$ and $A \subseteq B$, then $V(A) \leq V(B)$;
   (iii) If $A_n \uparrow A$, then $V(A_n) \uparrow V(A)$, $n \to \infty$.

4. Let $V$ be a capacity, $\xi \in L^2(\Omega, \mathcal{F}, P)$ and denote $C(\xi)$ by

$$C(\xi) := \int_{-\infty}^0 (V(\xi \geq t) - 1) \, dt + \int_0^\infty V(\xi \geq t) \, dt.$$

We call $C(\xi)$ the Choquet expectation of $\xi$ with respect to capacity $V$.

Dellacherie (1970) showed that comonotonic additivity is a necessary condition for a functional to be represented by a Choquet expectation.

*BSDEs and $g$-expectation.* Let $(\Omega, \mathcal{F}, P)$ be a probability space with filtration $(\mathcal{F}_s)_{s \geq 0}$, and let $(W_s)_{s \geq 0}$ be a standard $d$-Brownian motion. For ease of exposition, we assume $d = 1$. The results of this paper can be easily extended to the case $d > 1$. Suppose that $(\mathcal{F}_s)$ is the $\sigma$-filtration generated by $(W_s)_{s \geq 0}$, that is,

$$\mathcal{F}_s = \sigma\{W_r; 0 \leq r \leq s\}.$$

Let $T > 0$, $\mathcal{F}_T = \mathcal{F}$ and $g = g(y, z, t) : \mathbf{R} \times \mathbf{R}^d \times [0, T] \to \mathbf{R}$ be a function satisfying

(H.1) $\forall (y, z) \in \mathbf{R} \times \mathbf{R}^d$, $g(y, z, t)$ is continuous in $t$ and $\int_0^T g^2(0, 0, t)\, dt < \infty$;

(H.2) $g$ is uniformly Lipschitz continuous in $(y, z)$, that is, there exists a constant $C > 0$ such that $\forall y, y' \in \mathbf{R}$, $z, z' \in \mathbf{R}^d$, $|g(y, z, t) - g(y', z', t)| \leq C(|y - y'| + |z - z'|)$;

(H.3) $g(y, 0, t) \equiv 0$, $\forall (y, t) \in \mathbf{R} \times [0, T]$.

Let $\mathcal{M}(0, T, \mathbf{R}^n)$ be the set of all square integrable $\mathbf{R}^n$-valued, $\mathcal{F}_t$-adapted processes $\{v_t\}$ with

$$E \int_0^T |v_t|^2 \, dt < \infty.$$

For each $t \in [0, T]$, let $L^2(\Omega, \mathcal{F}_t, P)$ be the set of all $\mathcal{F}_t$-measurable random variables. Pardoux and Peng (1990) considered the following backward stochastic differential equation:

$$(1) \qquad y_t = \xi + \int_t^T g(y_s, z_s, s)\, ds - \int_t^T z_s\, dW_s,$$

and showed the following result:

LEMMA 1. *Suppose that $g$ satisfies* (H.1)–(H.3) *and $\xi \in L^2(\Omega, \mathcal{F}, P)$. Then BSDE* (1) *has a unique solution $(y, z) \in \mathcal{M}(0, T; \mathbf{R}) \times \mathcal{M}(0, T; \mathbf{R}^d)$.*

Using the solution of BSDE (1), Peng (1997) introduced the concept of $g$-expectation via BSDE (1).

DEFINITION 2. *Suppose $g$ satisfies* (H.1)–(H.3). *Given $\xi \in L^2(\Omega, \mathcal{F}, P)$, let $(y, z)$ be the solution of BSDE* (1). *We denote Peng's $g$-expectation of $\xi$ by $\mathcal{E}_g[\xi]$ and define it*

$$\mathcal{E}_g[\xi] := y_0.$$

From the definition of $g$-expectation, Peng (1997) introduced the concept of conditional $g$-expectation:

LEMMA 2.  *For any $\xi \in L^2(\Omega, \mathcal{F}, P)$, there exists unique $\eta \in L^2(\Omega, \mathcal{F}_t, P)$ such that*

$$\mathcal{E}_g[I_A \xi] = \mathcal{E}_g[I_A \eta] \qquad \forall A \in \mathcal{F}_t.$$

*We call $\eta$ the conditional $g$-expectation of $\xi$ and write $\eta$ as $\mathcal{E}_g[\xi|\mathcal{F}_t]$. Of course, such conditional expectations can be defined only for sub $\sigma$-algebras which appear in the filtration $\{\mathcal{F}_t\}$. Furthermore, $\mathcal{E}_g[\xi|\mathcal{F}_t]$ is the value of the solution $\{y_t\}$ of BSDE* (1) *at time $t$. That is,*

$$\mathcal{E}_g[\xi|\mathcal{F}_t] = y_t.$$

The $g$-expectation $\mathcal{E}_g[\cdot]$ preserves many of the properties of classical mathematical expectation. However, it does not preserve linearity. See, for example, Peng (1997) or Briand et al. (2000) for details.

LEMMA 3.  1. *If $c$ is a constant, then $\mathcal{E}_g[c] = c$.*
2. *If $\xi_1 \geq \xi_2$, then $\mathcal{E}_g[\xi_1] \geq \mathcal{E}_g[\xi_2]$.*
3. *$\mathcal{E}_g[\mathcal{E}_g[\xi|\mathcal{F}_t]] = \mathcal{E}_g[\xi]$.*
4. *If $\xi$ is $\mathcal{F}_t$-measurable, then $\mathcal{E}_g[\xi|\mathcal{F}_t] = \xi$.*
5. *For the real function $g$ defined on $\mathbf{R} \times \mathbf{R}^d \times [0, T]$, if $\xi$ is independent of $\mathcal{F}_t$, then $\mathcal{E}_g[\xi|\mathcal{F}_t] = \mathcal{E}_g[\xi]$.*

From the definition of $g$-expectation, it is natural to define $g$-probability thus:

DEFINITION 3.  For given $A \in \mathcal{F}$, denote $P_g(A)$ by

$$P_g(A) = \mathcal{E}_g[I_A].$$

We call $P_g(A)$ the $g$-probability of $A$. Obviously, $P_g(\cdot)$ is a capacity.

To simplify notation, we sometimes rewrite $g$-expectation $\mathcal{E}_g[\cdot]$, conditional $g$-expectation $\mathcal{E}_g[\cdot|\mathcal{F}_t]$ and $g$-probability $P_g(\cdot)$ as $\mathcal{E}_\mu[\cdot]$, $\mathcal{E}_\mu[\cdot|\mathcal{F}_t]$, $P_\mu(\cdot)$, respectively, if $g(y, z, t) = \mu_t|z|$.

REMARK 1.  1. $g$-expectation and conditional $g$-expectation depend on the choice of the function $g$, if $g$ is nonlinear, then $g$-expectation is usually also nonlinear.

2. If $g \equiv 0$, setting conditional expectation $E[\cdot|\mathcal{F}_t]$ on both sides of BSDE (1) yields $y_t = \mathcal{E}_g[\xi|\mathcal{F}_t] = E[\xi|\mathcal{F}_t]$, $y_0 = \mathcal{E}_g[\xi] = E[\xi]$. This implies another explanation for mathematical expectation: Within the particular framework of a Brownian filtration, conditional mathematical expectations with respect to $\mathcal{F}_t$ are the solution of a simple BSDE and mathematical expectation is the value of this solution at time $t = 0$.

The following example is a special case of Theorem 2.2 in Chen and Epstein ([2002](#)).

EXAMPLE 1.   Let $\mu := \{\mu_t\}$ be a continuous function on $[0, T]$. Suppose $g(y, z, t) = \mu_t |z|$ and $\mathcal{P}$ is a set of probability measures denoted by

$$\mathcal{P} := \left\{ Q^v : \frac{dQ^v}{dP} := e^{-(1/2) \int_0^T |v_s|^2 \, ds + \int_0^T v_s \, dW_s}, \right.$$

$$\left. v_t \text{ is } \mathcal{F}_t\text{-adapted and } |v_t| \le |\mu_t|, \text{ a.e. } t \in [0, T] \right\}.$$

Then for any $\xi \in L^2(\Omega, \mathcal{F}, P)$, we have:

(i)  $g$-expectation $\mathcal{E}_\mu[\xi]$:

$$\mathcal{E}_\mu[\xi] = \begin{cases} E\xi, & \mu = 0, \\ \inf_{Q \in \mathcal{P}} E_Q[\xi], & \mu < 0, \\ \sup_{Q \in \mathcal{P}} E_Q[\xi], & \mu > 0. \end{cases}$$

(ii)  Conditional $g$-expectation:

$$\mathcal{E}_\mu[\xi | \mathcal{F}_t] = \begin{cases} E[\xi | \mathcal{F}_t], & \mu = 0, \\ \operatorname*{ess\,inf}_{Q \in \mathcal{P}} E_Q[\xi | \mathcal{F}_t], & \mu < 0, \\ \operatorname*{ess\,sup}_{Q \in \mathcal{P}} E_Q[\xi | \mathcal{F}_t], & \mu > 0. \end{cases}$$

(iii)  $g$-probability (capacity): $\forall A \in \mathcal{F}$,

$$P_\mu(A) = \begin{cases} P(A), & \mu = 0, \\ \inf_{Q \in \mathcal{P}} Q(A), & \mu < 0, \\ \sup_{Q \in \mathcal{P}} Q(A), & \mu > 0. \end{cases}$$

The following are two key lemmas that we shall use in the next section: Lemma [4](#) is from Briand et al. ([2000](#)). We rewrite it in the following form. Lemma [5](#) is from Peng ([1997](#)):

LEMMA 4.   *Suppose* $\{X_t\}$ *is of the following form:*

$$dX_t = b_t \, dW_t,$$

$b$ *is a continuous, bounded adapted processes. Then*

$$\lim_{s \to t} \frac{\mathcal{E}_g[X_s | \mathcal{F}_t] - E[X_s | \mathcal{F}_t]}{s - t} = g(X_t, b_t, t),$$

*where the limit is in the sense of* $L^2(\Omega, \mathcal{F}_t, P)$.

LEMMA 5. *If $g$ is convex (resp. concave) in $(y, z)$, then for any $\xi, \eta \in L^2(\Omega, \mathcal{F}, P)$,*

$$\mathcal{E}_g[\xi + \eta | \mathcal{F}_t] \leq \quad (resp. \geq) \quad \mathcal{E}_g[\xi | \mathcal{F}_t] + \mathcal{E}_g[\eta | \mathcal{F}_t], \qquad t \in [0, T].$$

## 3. Main result.

The main result in this paper is the following theorem:

THEOREM 1. *Suppose that $g$ satisfies* (H.1)–(H.3). *Then there exists a Choquet expectation whose restriction to $L^2(\Omega, \mathcal{F}, P)$ is equal to a $g$-expectation if and only if $g$ does not depend on $y$ and is linear in $z$, that is, there exists a continuous function $\nu(t)$ such that*

$$g(y, z, t) = \nu(t) z.$$

The strategy of the proof is the following. First, we shall show that if $\mathcal{E}_g[\cdot]$ is a Choquet expectation on the set of all random variables with the form $y + z W_T$, then $g$ is of the form $g(y, z, t) = \mu_t |z| + \nu_t z$. Second, we further show if $g$-expectation is a Choquet expectation on the set of all random variables with the form $I_{[W_T \geq 1]}$ and $I_{[1 \leq W_T \leq 2]}$, then $\mu_t = 0$.

Lemma 6 is the first step. The first part of Lemma 6 shows the uniqueness of capacity:

LEMMA 6. *If there exists a Choquet capacity $V$ such that the associated Choquet expectation on $L^2(\Omega, \mathcal{F}, P)$ is equal to a $g$-expectation, then:*

(i) *$V(A) = P_g(A)$, $\forall A \in \mathcal{F}$;*

(ii) *There exist two continuous functions $\mu_t$, $\nu(t)$ on $[0, T]$ such that $g$ is of the form*

$$g(y, z, t) = \mu_t |z| + \nu(t) z.$$

PROOF. (i) Let $C(\xi)$ be the Choquet expectation of $\xi$ with respect to $V$, if

$$\mathcal{E}_g[\xi] = C(\xi) \qquad \forall \xi \in L^2(\Omega, \mathcal{F}, P).$$

In particular, for any $A \in \mathcal{F}$, let us choose $\xi = I_A$, thus, $\mathcal{E}_g[I_A] = C(I_A)$. By the definition of Choquet expectation, $C(I_A) = V(A)$ and $P_g(A) = \mathcal{E}_g[I_A]$, completing the proof of (i).

(ii) If $\mathcal{E}_g[\cdot]$ is a Choquet expectation on $L^2(\Omega, \mathcal{F}, P)$, then by Dellacherie's theorem in Dellacherie (1970). $\mathcal{E}_g[\cdot]$ is comonotonic additive. That is,

(2) $\quad \mathcal{E}_g[\xi + \eta] = \mathcal{E}_g[\xi] + \mathcal{E}_g[\eta] \qquad$ whenever $\xi$ and $\eta$ are comonotonic.

Choose constants $(y_1, z_1, t)$, $(y_2, z_2, t) \in \mathbf{R}^2 \times [0, T]$ such that $z_1 z_2 \geq 0$. For any $\tau \in [t, T]$, denote $\xi = y_1 + z_1(W_\tau - W_t)$ and $\eta = y_2 + z_2(W_\tau - W_t)$.

It is easy to check that $\xi$ and $\eta$ are comonotonic and independent of $\mathcal{F}_t$. Note that $g$ is deterministic and $y_i$ and $z_i$, $i = 1, 2$, are constants. Applying part 5 of Lemma 3,

$$\mathcal{E}_g[\xi|\mathcal{F}_t] = \mathcal{E}_g[\xi], \qquad \mathcal{E}_g[\eta|\mathcal{F}_t] = \mathcal{E}_g[\eta], \qquad \mathcal{E}_g[\xi + \eta|\mathcal{F}_t] = \mathcal{E}_g[\xi + \eta].$$

This with (2) implies

$$
\begin{aligned}
(3) \qquad & \frac{\mathcal{E}_g[\xi + \eta|\mathcal{F}_t] - E[\xi + \eta|\mathcal{F}_t]}{\tau - t} \\
& = \frac{\mathcal{E}_g[\xi|\mathcal{F}_t] - E[\xi|\mathcal{F}_t]}{\tau - t} + \frac{\mathcal{E}_g[\eta|\mathcal{F}_t] - E[\eta|\mathcal{F}_t]}{\tau - t}.
\end{aligned}
$$

Let $\tau \to t$ on both sides of (3) to obtain, using Lemma 4,

$$
\begin{aligned}
(4) \qquad g(y_1 + y_2, z_1 + z_2, t) &= g(y_1, z_1, t) + g(y_2, z_2, t) \\
& \forall z_1 z_2 \geq 0, y_1, y_2 \in \mathbf{R},
\end{aligned}
$$

which implies that $g$ is linear with respect to $y$ in $\mathbf{R}$ and $z$ in $\mathbf{R}_+$ (or $\mathbf{R}_-$).

Thus, for any $(y, z, t) \in \mathbf{R}^2 \times [0, T]$, note that $g(y, 0, t) = 0$ in (H.3) and apply (4) to obtain

$$
\begin{aligned}
g(y, z, t) &= g(y + 0, zI_{[z \geq 0]} + zI_{[z \leq 0]}, t) \\
&= g(y, zI_{[z \geq 0]}, t) + g(0, zI_{[z \leq 0]}, t) \\
&= g(y + 0, 0 + zI_{[z \geq 0]}, t) + g(0, -(-z)I_{[z \leq 0]}, t) \\
&= g(y, 0, t) + g(0, zI_{[z \geq 0]}, t) + g(0, -(-z)I_{[z \leq 0]}, t) \\
&= g(0, 1, t)zI_{[z \geq 0]} - g(0, -1, t)zI_{[z \leq 0]} \\
&= g(0, 1, t)z^+ + g(0, -1, t)(-z)^+ \\
&= g(0, 1, t)\frac{|z| + z}{2} + g(0, -1, t)\frac{|z| - z}{2} \\
&= \frac{g(0, 1, t) + g(0, -1, t)}{2}|z| + \frac{g(0, 1, t) - g(0, -1, t)}{2}z.
\end{aligned}
$$

The second equality is because $zI_{[z \geq 0]} \cdot zI_{[z \leq 0]} = 0$.

Set $\mu_t := \frac{g(0,1,t) + g(0,-1,t)}{2}$ and $\nu(t) := \frac{g(0,1,t) - g(0,-1,t)}{2}$ to complete the proof. $\qquad \square$

Next, we show that $\mu_t = 0$ for $t \in [0, T]$. We need the following lemmas. Lemma 7 is a special case of the comonotonic theorem in Chen, Kulperger and Wei (2005):

LEMMA 7. *Suppose* $\Phi$ *is a function such that* $\Phi(W_T) \in L^2(\Omega, \mathcal{F}, P)$. *Let* $(y_t, z_t)$ *be the solution of*

$$(5) \qquad y_t = \Phi(W_T) + \int_t^T \mu_t |z_s| \, ds - \int_t^T z_s \, dW_s,$$

*where* $\mu_t$ *is a continuous function on* $[0, T]$.

  (i) *If* $\Phi$ *is increasing, then* $z_t \geq 0$, *a.e.* $t \in [0, T]$.
  (ii) *If* $\Phi$ *is decreasing, then* $z_t \leq 0$, *a.e.* $t \in [0, T]$.

PROOF (Sketched for the reader's convenience). For $\varepsilon > 0$, let $g_\varepsilon(z, t) = \mu_t \sqrt{z^2 + \varepsilon}$, then $g_\varepsilon$ is a smooth $C^2$-function which $\to \mu_t |z|$ as $\varepsilon \to 0$.

Let $\{y_s^{\varepsilon, t, x}, z_s^{\varepsilon, t, x}\}_{(0 \leq s \leq T)}$ be the solution of the BSDE:

$$y_s = \Phi(W_T - W_t + x) + \int_s^T \mu_r \sqrt{z_r 2 + \varepsilon} \, dr - \int_s^T z_r \, dW_r, \qquad 0 \leq s \leq T,$$

and $\{y_s^{t, x}, z_s^{t, x}\}_{(0 \leq s \leq T)}$ be the solution of the BSDE:

$$y_s = \Phi(W_T - W_t + x) + \int_s^T \mu_r |z_r| \, dr - \int_s^T z_r \, dW_r, \qquad 0 \leq s \leq T.$$

By the convergence theorem of BSDE, see Proposition 2.1 in El Karoui, Peng and Quenez ([1997](#)),

$$\{y_s^{\varepsilon, t, x}, z_s^{\varepsilon, t, x}\}_{(0 \leq s \leq T)} \to \{y_s^{t, x}, z_s^{t, x}\}_{(0 \leq s \leq T)} \qquad \text{as } \varepsilon \to 0$$

in the sense of $\mathcal{M}(0, T; \mathbf{R}) \times \mathcal{M}(0, T; \mathbf{R}^d)$.

Moreover, if we choose $x = 0$, $t = 0$ in $\{y_s^{t, x}, z_s^{t, x}\}$, then $\{y_s^{0, 0}, z_s^{0, 0}\}$ is the solution of BSDE ([5](#)). Thus, if we can show for each $t \in [0, T]$, $z_s^{\varepsilon, t, x} \geq 0$, a.e. $s \in [0, T]$, by the convergence theorem of BSDE in El Karoui, Peng and Quenez ([1997](#)), we have $z_s^{t, x} \geq 0$, thus, $z_t = z_s^{0, 0} \geq 0$.

Without loss of generality, we assume that $\Phi$ is a smooth $C^2$-function, otherwise, we can choose a sequence of smooth $C^2$-functions $\Phi_\varepsilon$ such that $\Phi_\varepsilon \to \Phi$ as $\varepsilon \to 0$.

Let $u_\varepsilon(t, x) = y_t^{\varepsilon, t, x}$. Then the general Feynman–Kac formula, see Proposition 4.3 in El Karoui, Peng and Quenez ([1997](#)) or Ma, Protter and Yong ([1994](#)), implies that $u_\varepsilon$ solves the PDE

$$\frac{\partial u_\varepsilon}{\partial t} + \frac{1}{2} \frac{\partial^2 u_\varepsilon}{\partial x^2} + g_\varepsilon(x, t) = 0,$$

$$u_\varepsilon(T, x) = \Phi(x), 0 \leq t \leq T.$$

Moreover,

$$z_s^{\varepsilon, t, x} = \frac{\partial u_\varepsilon(s, W_s - W_t)}{\partial x}.$$

On the other hand, by the comparison theorem of PDE, $u_\varepsilon(t,x)$ is increasing in $x$ if $\Phi$ is increasing, thus,

$$\frac{\partial u_\varepsilon(t,x)}{\partial x} \geq 0.$$

This implies that $z_s^{\varepsilon,t,x} \geq 0, s \geq 0$. Letting $\varepsilon \to 0$ and $(x,t) = (0,0)$, we obtain (i). Similarly, we can obtain (ii) if $\Phi$ is decreasing. The proof is complete. □

Furthermore, we can prove the following:

LEMMA 8. *Let $\mu_t$ be a continuous function on $[0,T]$ and $(y,z)$ be the solution of BSDE*

$$(6) \qquad\qquad y_t = \xi + \int_t^T \mu_s |z_s| \, ds - \int_t^T z_s \, dW_s.$$

   (i) *If $\xi = I_{[W_T \geq 1]}$, then $z_t > 0 \ \forall t \in [0,T)$.*
   (ii) *If $\xi = \Phi(W_T)$, where $\Phi$ is a bounded function with strictly positive derivative $\Phi'$, then $z_t > 0, \ \forall t \in [0,T)$.*
   (iii) *If $\xi = I_{[2 \geq W_T \geq 1]}$, then $P \times \lambda(\{(\omega,t) : z_t(\omega) < 0\}) > 0$,*
*where $\lambda$ denotes Lebesgue measure on $[0,T)$ and $P \times \lambda$ denotes the product of the probability measure $P$ and the Lebesgue measure $\lambda$.*

PROOF.   (i) Note that the indicator function $I_{[x \geq 1]}$ is increasing, so by Lemma 7, $z_t \geq 0$, a.e. $t \in [0,T]$. This implies that the BSDE (6) is actually a linear BSDE:

$$y_t = I_{[W_T \geq 1]} + \int_t^T \mu_s z_s \, ds - \int_t^T z_s \, dW_s.$$

Let

$$\overline{W}_t = W_t - \int_0^t \mu_s \, ds,$$

then

$$(7) \qquad\qquad y_t = I_{[W_T \geq 1]} - \int_t^T z_s \, d\overline{W}_s.$$

Let $Q$ be the probability measure defined by

$$\frac{dQ}{dP} = \exp\left[-\frac{1}{2}\int_0^T \mu_s^2 \, ds + \int_0^T \mu_s \, dW_s\right].$$

By Girsanov's lemma, $(\overline{W}_t)_{0 \leq t \leq T}$ is a $Q$-Brownian motion.

Set conditional expectation $E_Q[\cdot|\mathcal{F}_t]$ on both sides of BSDE (7). From the Markov property,

$$
\begin{aligned}
y_t &= E_Q[I_{[W_T \geq 1]}|\mathcal{F}_t] \\
&= E_Q[I_{[\overline{W}_T \geq 1 - \int_0^T \mu_s\,ds]}|\mathcal{F}_t] \\
&= E_Q[I_{[\overline{W}_T - \overline{W}_t \geq 1 - \int_0^T \mu_s\,ds - \overline{W}_t]}|\mathcal{F}_t] \\
&= E_Q[I_{[\overline{W}_T - \overline{W}_t \geq 1 - \int_0^T \mu_s\,ds - \overline{W}_t]}|\sigma(W_t)].
\end{aligned}
$$

Note that $\sigma(W_s; s \leq t) = \sigma(\overline{W}_s; s \leq t)$ because $\mu_t$ is a real function in $t$, thus,

$$
y_t = E_Q[I_{[\overline{W}_T - \overline{W}_t \geq 1 - \int_0^T \mu_s\,ds - \overline{W}_t]}|\sigma(\overline{W}_t)].
$$

Since $\overline{W}_T - \overline{W}_t$ and $\overline{W}_t$ are independent, we have

$$
y_t = E_Q[I_{[\overline{W}_T - \overline{W}_t \geq 1 - \int_0^T \mu_s\,ds - h]}]|_{h = \overline{W}_t}.
$$

But $\overline{W}_T - \overline{W}_t \sim N(0, T - t)$, therefore,

$$
y_t = \int_{1 - \int_0^T \mu_s\,ds - h}^{\infty} \varphi(x)\,dx\Big|_{h = \overline{W}_t},
$$

where $\varphi(x)$ is the density function of the normal distribution $N(0, T - t)$.

By the relation between $y_t$ and $z_t$, see Corollary 4.1 in El Karoui, Peng and Quenez (1997), we have

$$
z_t = \frac{\partial y_t}{\partial h}\Big|_{h = \overline{W}_t} = \varphi\left(1 - \int_0^T \mu_s\,ds - \overline{W}_t\right) > 0,
$$

that is, $z_t > 0$, a.e. $t \in [0, T)$.

(ii) Similar to the proof of (i) and noting the fact that $\Phi' > 0$, it is easy to check

$$
z_t = \int_{-\infty}^{\infty} \Phi'\left(x + \overline{W}_t + \int_0^T \mu_s\,ds\right)\varphi(x)\,dx > 0, \qquad t \in [0, T),
$$

where $\varphi(x)$ is the density function of the normal distribution $N(0, T - t)$.

(iii) For given $\xi = I_{[2 \geq W_T \geq 1]}$, we assume the conclusion of (ii) is false, then $z_t \geq 0$, a.e. $t \in [0, T)$, which implies that BSDE (6) is actually a linear BSDE:

$$
y_t = I_{[2 \geq W_T \geq 1]} + \int_t^T \mu_s z_s\,ds - \int_t^T z_s\,dW_s.
$$

That is,

$$
(8) \qquad y_t = I_{[2 \geq W_T \geq 1]} - \int_t^T z_s\,d\overline{W}_s,
$$

where

$$\overline{W}_t = W_t - \int_0^t \mu_s \, ds.$$

As in (i), let

$$\frac{dQ}{dP} = \exp\left[-\frac{1}{2}\int_0^T \mu_s^2 \, ds + \int_0^T \mu_s \, dW_s\right].$$

Applying Girsanov's lemma again, $(\overline{W}_t)_{0 \le t \le T}$ is a $Q$-Brownian motion.

Set conditional expectation $E_Q[\cdot|\mathcal{F}_t]$ on both sides of BSDE (8). Note that $\sigma(W_s; s \le t) = \sigma(\overline{W}_s; s \le t)$,

$$\begin{aligned}
y_t &= E_Q[I_{[2 \ge W_T \ge 1]}|\mathcal{F}_t] \\
&= E_Q[I_{[2 - \int_0^T \mu_s \, ds - \overline{W}_T \ge \overline{W}_T - \overline{W}_t \ge 1 - \int_0^T \mu_s \, ds - \overline{W}_t]}|\mathcal{F}_t] \\
&= E_Q[I_{[2 - \int_0^T \mu_s \, ds - \overline{W}_T \ge \overline{W}_T - \overline{W}_t \ge 1 - \int_0^T \mu_s \, ds - \overline{W}_t]}|\sigma(\overline{W}_t)] \\
&= E_Q[I_{[2 - \int_0^T \mu_s \, ds - h \ge \overline{W}_T - \overline{W}_t \ge 1 - \int_0^T \mu_s \, ds - h]}]|_{h = \overline{W}_t}.
\end{aligned}$$

Since $\overline{W}_T - \overline{W}_t \sim N(0, T - t)$,

$$y_t = \int_{1 - \int_0^T \mu_s \, ds - h}^{2 - \int_0^T \mu_s \, ds - h} \varphi(x) \, dx \Big|_{h = \overline{W}_t}.$$

Therefore, applying the relation between $y_t$ and $z_t$ again,

$$\begin{aligned}
z_t &= \frac{\partial y_t}{\partial h}\Big|_{h = \overline{W}_t} \\
&= \varphi\left(1 - \int_0^T \mu_s \, ds - \overline{W}_t\right) - \varphi\left(2 - \int_0^T \mu_s \, ds - \overline{W}_t\right) \\
&= \frac{1}{\sqrt{2\pi(T - t)}} \exp\left[-\frac{(1 - \int_0^T \mu_s \, ds - \overline{W}_t)^2}{2(T - t)}\right] \\
&\quad - \frac{1}{\sqrt{2\pi(T - t)}} \exp\left[-\frac{(2 - \int_0^T \mu_s \, ds - \overline{W}_t)^2}{2(T - t)}\right].
\end{aligned}$$

However, it is easy to check that

$$z_t > 0, t \in [0, T) \qquad \text{when } \overline{W}_t < \frac{3}{2} - \int_0^T \mu_s \, ds;$$

$$z_t < 0, t \in [0, T) \qquad \text{when } \overline{W}_t > \frac{3}{2} - \int_0^T \mu_s \, ds,$$

which implies

$$P(z_t > 0) > 0, \qquad P(z_t < 0) > 0 \qquad \text{a.e. } t \in [0, T),$$

thus, $P \times \lambda((\omega, t) : z_t(\omega) < 0) > 0$. We obtain a contradiction. The proof is complete. $\square$

Let $L_+^2(\Omega, \mathcal{F}, P)$ [resp. $L_-^2(\Omega, \mathcal{F}, P)$] be the set of all nonnegative (resp. nonpositive) random variables in $L^2(\Omega, \mathcal{F}, P)$.

LEMMA 9. *Suppose that $g$ is a convex (or concave) function. If $\mathcal{E}_g[\cdot]$ is comonotonic additive on $L_+^2(\Omega, \mathcal{F}, P)$ [resp. $L_-^2(\Omega, \mathcal{F}, P)$], then $\mathcal{E}_g[\cdot|\mathcal{F}_t]$ is also comonotonic additive on $L_+^2(\Omega, \mathcal{F}, P)$ [resp. $L_-^2(\Omega, \mathcal{F}, P)$] for any $t \in [0, T)$.*

PROOF. We show the above result on $L_+^2(\Omega, \mathcal{F}, P)$; the result on $L_-^2(\Omega, \mathcal{F}, P)$ can be proved in the same way.

Because $\mathcal{E}_g[\cdot]$ is comonotonic additive on $L_+^2(\Omega, \mathcal{F}, P)$, then $\forall \xi, \eta \in L_+^2(\Omega, \mathcal{F}, P)$, we have

(9)     $\mathcal{E}_g[\xi + \eta] = \mathcal{E}_g[\xi] + \mathcal{E}_g[\eta]$ whenever $\xi$ and $\eta$ are comonotonic.

We now show $\forall t \in [0, T)$

(10)   $\mathcal{E}_g[\xi + \eta|\mathcal{F}_t] = \mathcal{E}_g[\xi|\mathcal{F}_t] + \mathcal{E}_g[\eta|\mathcal{F}_t]$ whenever $\xi$ and $\eta$ are comonotonic.

First, we assume that $g$ is a convex function. Then by Lemma 5,

$$\mathcal{E}_g[\xi + \eta|\mathcal{F}_t] \le \mathcal{E}_g[\xi|\mathcal{F}_t] + \mathcal{E}_g[\eta|\mathcal{F}_t] \qquad \forall t \in [0, T].$$

If (10) is false, then there exists $t \in [0, T]$ such that

$$P(\omega : \mathcal{E}_g[\xi + \eta|\mathcal{F}_t] < \mathcal{E}_g[\xi|\mathcal{F}_t] + \mathcal{E}_g[\eta|\mathcal{F}_t]) > 0.$$

Let

$$A = \{\omega : \mathcal{E}_g[\xi + \eta|\mathcal{F}_t] < \mathcal{E}_g[\xi|\mathcal{F}_t] + \mathcal{E}_g[\eta|\mathcal{F}_t]\}.$$

Obviously $A \in \mathcal{F}$, and

$$I_A \mathcal{E}_g[\xi + \eta|\mathcal{F}_t] < I_A \mathcal{E}_g[\xi|\mathcal{F}_t] + I_A \mathcal{E}_g[\eta|\mathcal{F}_t].$$

Set $g$-expectation $\mathcal{E}_g[\cdot]$ on both sides of the above inequality. By the strict comparison theorem of BSDE, see Peng (1997), we have

(11)   $\mathcal{E}_g[I_A \mathcal{E}_g[\xi + \eta|\mathcal{F}_t]] < \mathcal{E}_g\{I_A \mathcal{E}_g[\xi|\mathcal{F}_t] + I_A \mathcal{E}_g[\eta|\mathcal{F}_t]\}.$

Observing the above inequality, apply the convexity of $g$ again to the right-hand side of (11) and use part 3 of Lemma 3 to obtain

$$\mathcal{E}_g\{I_A \mathcal{E}_g[\xi|\mathcal{F}_t] + I_A \mathcal{E}_g[\eta|\mathcal{F}_t]\} \le \mathcal{E}_g[I_A \mathcal{E}_g[\xi|\mathcal{F}_t]] + \mathcal{E}_g[I_A \mathcal{E}_g[\eta|\mathcal{F}_t]]$$
$$= \mathcal{E}_g[I_A \xi] + \mathcal{E}_g[I_A \eta].$$

But the left-hand side of (11) is

$$\mathcal{E}_g[I_A\mathcal{E}_g[\xi + \eta|\mathcal{F}_t]] = \mathcal{E}_g[I_A\xi + I_A\eta].$$

Thus,

$$(12) \qquad \mathcal{E}_g[I_A\xi + I_A\eta] < \mathcal{E}_g[I_A\xi] + \mathcal{E}_g[I_A\eta].$$

Furthermore, since $\xi$ and $\eta$ are positive and comonotonic, obviously $I_A\xi$ and $I_A\eta$ are also positive and comonotonic, by the assumption that $\mathcal{E}_g[\cdot]$ is comonotonic additive, and Dellacherie's (1970) theorem,

$$(13) \qquad \mathcal{E}_g[I_A\xi + I_A\eta] = \mathcal{E}_g[I_A\xi] + \mathcal{E}_g[I_A\eta].$$

Inequality (12) contradicts (13), thus,

$$\mathcal{E}_g[\xi + \eta|\mathcal{F}_t] = \mathcal{E}_g[\xi|\mathcal{F}_t] + \mathcal{E}_g[\eta|\mathcal{F}_t] \qquad \forall\, t \in [0, T].$$

Next, if $g$ is concave, then, by Lemma 5,

$$\mathcal{E}_g[\xi + \eta|\mathcal{F}_t] \geq \mathcal{E}_g[\xi|\mathcal{F}_t] + \mathcal{E}_g[\eta|\mathcal{F}_t] \qquad \forall\, t \in [0, T].$$

The rest can be proved in a fashion similar to result (i).  □

Combining Dellacherie's (1970) theorem and Lemma 9, we immediately obtain the following:

COROLLARY 1.  *Under the assumption of Lemma 9, if $\mathcal{E}_g[\cdot]$ is a Choquet expectation on $L^2_+(\Omega, \mathcal{F}, P)$ [resp. $L^2_-(\Omega, \mathcal{F}, P)$], then for each $t \in [0, T]$, $\mathcal{E}_g[\cdot|\mathcal{F}_t]$ is also a Choquet expectation on $L^2_+(\Omega, \mathcal{F}, P)$ [resp. $L^2_-(\Omega, \mathcal{F}, P)$].*

We now study the case where $g$ is of the form $g(t, y, z) = \mu_t|z|$. Obviously, if $\mu_t \geq 0$, $t \in [0, T]$, then $g$ is convex and if $\mu_t \leq 0$, $t \in [0, T]$, $g$ is concave.

LEMMA 10.  *Let $\mu_t \neq 0$ be a continuous function on $[0, T]$ and $g(z, t) = \mu_t|z|$. There exists no Choquet expectation agreeing with $\mathcal{E}_\mu[\cdot]$ on $L^2(\Omega, \mathcal{F}, P)$.*

PROOF.  Assume the result of this lemma is false, then $\mathcal{E}_\mu[\cdot]$ is a Choquet expectation on $L^2(\Omega, \mathcal{F}, P)$. By Dellacherie's (1970) theorem, $\mathcal{E}_\mu[\cdot]$ is comonotonic additive on $L^2(\Omega, \mathcal{F}, P)$.

We now choose two special random variables $\xi_1 = I_{[W_T \geq 1]}$ and $\xi_2 = I_{[2 \geq W_T \geq 1]}$. Let $(y^i, z^i)$, $i = 1, 2$, be the solutions of the following BSDEs corresponding to $\xi_1$ and $\xi_2$, respectively,

$$y_t = \xi_i + \int_t^T \mu_s|z_s|\,ds - \int_t^T z_s\,dW_s, \qquad i = 1, 2.$$

If $(\overline{y}_t, \overline{z}_t)$ is the solution of BSDE,

$$\overline{y}_t = (\xi_1 + \xi_2) + \int_t^T \mu_s |\overline{z}_s|\, ds - \int_t^T \overline{z}_s\, dW_s,$$

then $y_t^1 = \mathcal{E}_\mu[\xi_1 | \mathcal{F}_t]$, $y_t^2 = \mathcal{E}_\mu[\xi_2 | \mathcal{F}_t]$ and $\overline{y}_t = \mathcal{E}_\mu[\xi_1 + \xi_2 | \mathcal{F}_t]$.

It is easy to show that $\xi_1$ and $\xi_2$ are positive and comonotonic. As we have assumed that $g$-expectation $\mathcal{E}_\mu[\cdot]$ is a Choquet expectation, thus, $\mathcal{E}_\mu[\cdot]$ is comonotonic additive with respect to $\xi_1$, $\xi_2$. By Lemma 9, $\mathcal{E}_\mu[\cdot | \mathcal{F}_t]$ is also comonotonic additive with respect to $\xi_1$, $\xi_2$, that is,

$$\mathcal{E}_\mu[\xi_1 + \xi_2 | \mathcal{F}_t] = \mathcal{E}_\mu[\xi_1 | \mathcal{F}_t] + \mathcal{E}_\mu[\xi_2 | \mathcal{F}_t] \qquad \forall t \in [0, T].$$

This can be written in another form, namely

$$(14) \qquad \overline{y}_t = y_t^1 + y_t^2 \qquad \forall t \in [0, T].$$

Let $\langle X, W \rangle$ be the finite variation process generated by the semi-martingale $X$ and Brownian motion $W$, then, from (14),

$$\langle \overline{y}, W \rangle_t = \langle y^1 + y^2, W \rangle_t = \langle y^1, W \rangle_t + \langle y^2, W \rangle_t \qquad \forall t \in [0, T],$$

but

$$\overline{z}_t = \frac{d\langle \overline{y}, W \rangle_t}{dt}, \qquad z_t^1 = \frac{d\langle y^1, W \rangle_t}{dt}, \qquad z_t^2 = \frac{d\langle y^2, W \rangle_t}{dt}.$$

Thus,

$$\overline{z}_t = z_t^1 + z_t^2 \qquad \text{a.e. } t \in [0, T].$$

Applying the above inequality to (14), note that (14) can be rewritten as

$$(\xi_1 + \xi_2) + \int_t^T \mu_s |\overline{z}_s|\, ds - \int_t^T \overline{z}_s\, dW_s = \sum_{i=1}^2 \left( \xi_i + \int_t^T \mu_s |z_s^i|\, ds - \int_t^T z_s^i\, dW_s \right).$$

We can obtain

$$\mu_t |z_t^1 + z_t^2| = \mu_t |z_t^1| + \mu_t |z_t^2| \qquad \text{a.e. } t \in [0, T].$$

Since $\mu_t \neq 0$, therefore,

$$(15) \qquad |z_t^1 + z_t^2| = |z_t^1| + |z_t^2| \qquad \text{a.e. } t \in [0, T].$$

Obviously, (15) is true if and only if $z_t^1 z_t^2 \geq 0$.

However, from Lemma 8, we know $z_t^1 > 0$, a.e. $t \in [0, T]$ and

$$P \times \lambda((\omega, t) : z_t^2(\omega) < 0) > 0.$$

Thus, $P \times \lambda((\omega, t) : z_t^1(\omega) z_t^2(\omega) < 0) > 0$, which implies

$$P \times \lambda((\omega, t) : |z_t^1(\omega) + z_t^2(\omega)| < |z_t^1(\omega)| + |z_t^2(\omega)|) > 0,$$

which contradicts (15). The lemma's proof is complete. $\square$

From the above proof, applying the strict comparison theorem of BSDE in Peng (1997), we have the following:

COROLLARY 2.  *If $\mu_t \neq 0$, $\forall t \in [0, T]$. Let $\xi_1 = I_{[W_T \geq 1]}$ and $\xi_2 = I_{[2 \geq W_T \geq 1]}$, obviously $\xi_1$ and $\xi_2$ are comonotonic, but $\mathcal{E}_\mu[\xi_1 + \xi_2] < \mathcal{E}_\mu[\xi_1] + \mathcal{E}_\mu[\xi_2]$.*

We now prove our main theorem.

PROOF OF THEOREM 1.

*Sufficiency.*  If $g(y, z, t) = \nu_t z$, for any $\xi \in L^2(\Omega, \mathcal{F}, P)$, let us consider BSDE

$$y_t = \xi + \int_t^T \nu_s z_s \, ds - \int_t^T z_s \, dW_s.$$

Let $\overline{W}_t = W_t - \int_0^t \nu_s \, ds$, then

$$y_t = \xi - \int_t^T z_s \, d\overline{W}_s.$$

By Girsanov's lemma, $(\overline{W}_t)_{0 \leq t \leq T}$ is a $Q$-Brownian motion under $Q$ defined by

$$\frac{dQ}{dP} = \exp\left[ -\frac{1}{2} \int_0^T v_s^2 \, ds + \int_0^T v_s \, dW_s \right].$$

Thus,

$$\mathcal{E}_g[\xi | \mathcal{F}_t] = E_Q[\xi | \mathcal{F}_t], \qquad \mathcal{E}_g[\xi] = E_Q[\xi].$$

This implies $g$-expectation is a classical mathematical expectation. Obviously, the classical mathematical expectation can be represented by the Choquet expectation. So the sufficiency proof is complete.

*Necessity.*  For any $\xi \in L^2(\Omega, \mathcal{F}, P)$, by Lemma 6(ii), there exist two continuous functions on $[0, T]$ such that

$$g(y, z, t) = \mu_t |z| + \nu(t) z.$$

Without loss of generality, we assume $\nu(t) = 0, t \in [0, T]$, otherwise, by Girsanov's lemma, we can rewrite the BSDE

$$y_t = \xi + \int_t^T (\mu_s |z_s| + \nu_s z_s) \, ds - \int_t^T z_s \, dW_s$$

as

(16) $$y_t = \xi + \int_t^T \mu_s |z_s| \, ds - \int_t^T z_s \, d\overline{W}_s,$$

where $\overline{W}_s := W_s - \int_0^s \nu(r) \, dr$, $(\overline{W}_t)_{0 \leq t \leq T}$ is a $Q$-Brownian motion under $Q$ defined by

$$\frac{dQ}{dP} = \exp\left[ -\frac{1}{2} \int_0^T \nu_s^2 \, ds + \int_0^T \nu_s \, dW_s \right].$$

We can consider our question on the probability space $(\Omega, \mathcal{F}, Q)$.

Assume $\mu \not\equiv 0$, then there exists $t_0$ such that $\mu_{t_0} \neq 0$. Without loss of generality, we assume $\mu_{t_0} > 0$.

Since $\mu_t$ is continuous, then there exists a region of $t_0$, say $[\overline{t}, \overline{T}] \subset [0, T]$ such that $\forall\, t \in [\overline{t}, \overline{T}]$, $\mu_t > 0$.

The next step of the proof is to localize in time so as to use Lemma 10.

Let $\xi_1 = I_{[W_{\overline{T}} - W_{\overline{t}} \geq 1]}$ and $\xi_2 = I_{[2 \geq W_{\overline{T}} - W_{\overline{t}} \geq 1]}$. Obviously, $\xi_1$ and $\xi_2$ are comonotonic.

We now show that

$$\mathcal{E}_\mu[\xi_1 + \xi_2] < \mathcal{E}_\mu[\xi_1] + \mathcal{E}_\mu[\xi_2],$$

which implies that $\mathcal{E}_\mu[\cdot]$ is not comonotonic additive for comonotonic random variables $\xi$ and $\eta$.

Let $\overline{W}_s = \overline{W}_{\overline{t}+s} - \overline{W}_{\overline{t}}$, then $\{\overline{W}_s : 0 \leq s \leq \overline{T} - \overline{t}\}$ is $(\mathcal{F}'_s)$ Brownian motion, where

$$\mathcal{F}'_s = \sigma\{\overline{W}_r : 0 \leq r \leq s\} = \sigma\{W_{\overline{t}+r} - W_{\overline{t}} : 0 \leq r \leq s\}.$$

Using the above notation, $\xi_1$ and $\xi_2$ can be rewritten as $\xi_1 = I_{[\overline{W}_{\overline{T}-\overline{t}} \geq 1]}$ and $\xi_2 = I_{[2 \geq \overline{W}_{\overline{T}-\overline{t}} \geq 1]}$.

For the given $\xi_1$ and $\xi_2$, let $a_t = \mu_{t+\overline{t}}$ and $(Y^i, Z^i)$ be the solutions of the following BSDEs with terminal value $\xi_1$ and $\xi_2$, respectively, on $[0, \overline{T} - \overline{t}]$:

$$(17) \qquad Y^i_t = \xi_i + \int_t^{\overline{T}-\overline{t}} a_s |Z^i_s|\, ds - \int_t^{\overline{T}-\overline{t}} Z^i_s\, d\overline{W}_s, \qquad t \in [0, \overline{T} - \overline{t}], i = 1, 2,$$

and $(\overline{Y}, \overline{Z})$ be the solution of the BSDE:

$$(18) \qquad \overline{Y}_t = \xi_1 + \xi_2 + \int_t^{\overline{T}-\overline{t}} a_s |\overline{Z}^i_s|\, ds - \int_t^{\overline{T}-\overline{t}} \overline{Z}^i_s\, d\overline{W}_s, \qquad t \in [0, \overline{T} - \overline{t}].$$

Since $a_t = \mu_{\overline{t}+t} \neq 0$, $\forall\, t \in [0, \overline{T} - \overline{t}]$, by Corollary 2,

$$(19) \qquad \qquad \overline{Y}_t < \overline{Y}^1_t + \overline{Y}^2_t, \qquad t \in [0, \overline{T} - \overline{t}].$$

On the other hand, for the given $\xi_1$ and $\xi_2$, consider the BSDE on $[0, T]$:

$$(20) \qquad y^i_t = \xi_i + \int_t^T \mu_s |z^i_s|\, ds - \int_t^T z^i_s\, dW_s, \qquad i = 1, 2, t \in [0, T],$$

and

$$(21) \qquad \overline{y}_t = \xi_1 + \xi_2 + \int_t^T \mu_s |\overline{z}^i_s|\, ds - \int_t^T \overline{z}^i_s\, dW_s, \qquad t \in [0, T].$$

Comparing (17) with (20) and (18) with (21),

$$Y^i_t = y^i_t, \qquad i = 1, 2; \qquad \overline{Y}_t = \overline{y}_t, \qquad t \in [0, \overline{T} - \overline{t}].$$

But $y_t^1 = \mathcal{E}_\mu[\xi_1|\mathcal{F}_t]$, $y_t^2 = \mathcal{E}_\mu[\xi_2|\mathcal{F}_t]$ and $\overline{y}_t = \mathcal{E}_\mu[\xi_1 + \xi_2|\mathcal{F}_t]$.

Thus,

$$Y_0^i = \mathcal{E}_\mu[\xi_i], \qquad i = 1, 2, \qquad \overline{Y}_0 = \mathcal{E}_\mu[\xi_1 + \xi_2].$$

Applying (19),

$$\mathcal{E}_\mu[\xi_1 + \xi_2] < \mathcal{E}_\mu[\xi_1] + \mathcal{E}_\mu[\xi_2],$$

which contradicts the comonotonic additivity of $\mathcal{E}_\mu[\cdot]$. Thus, $\mu(t) = 0$ $\forall t \in [0, T]$. The proof is complete. $\square$

An interesting application of Theorem 1 is the following:

COROLLARY 3. *Suppose $\mu \neq 0$ and let $\mathcal{E}_\mu[\cdot]$ be the maximal (minimal) expectations defined in Example 1, then maximal (minimal) expectation is not a Choquet expectation on $L^2(\Omega, \mathcal{F}, P)$.*

REMARK 2. 1. In the proofs of Lemma 6(ii) and Theorem 1, we only use random variables having the form $y + zW_t$ and $I_{[W_T \in (a,b)]}$. Thus, Lemma 6(ii) and Theorem 1 actually imply that if and only if $g$ is linear in $z$, then $g$-expectation is a Choquet expectation on the set of all random variables with the form $f(W_T) \in L^2(\Omega, \mathcal{F}, P)$.

2. Because $g$-expectation depends on the choice of $g$, if $g$ is nonlinear in $z$, Theorem 1 implies that $g$-expectation is not a Choquet expectation on $L^2(\Omega, \mathcal{F}, P)$.

3. It is well understood that mathematical expectation is linear in the sense of

$$E(\xi + \eta) = E\xi + E\eta \qquad \forall \xi, \eta \in L^2(\Omega, \mathcal{F}, P).$$

For the Choquet expectation, the above equality is still true when $\xi$ and $\eta$ are comonotonic. However, for the $g$-expectation, if $g$ is nonlinear, the above additivity no longer holds even for comonotonic random variables. From this viewpoint, our result implies Peng's $g$-expectation usually is more nonlinear than the Choquet expectation on $L^2(\Omega, \mathcal{F}, P)$.

## 4. Feynman–Kac formula and Choquet expectation. Let $u$ be the solution of the partial differential equation (PDE)

$$(22) \qquad \begin{aligned} \frac{\partial u(t,x)}{\partial t} &= \frac{1}{2}\frac{\partial^2 u(t,x)}{\partial x^2}, \\ u(0,x) &= f(x), \qquad t \geq 0, x \in R. \end{aligned}$$

By the famous Feynman–Kac formula, the solution $u(t,x)$ of PDE (22) can be represented by mathematical expectation:

$$(23) \qquad u(t,x) = Ef(W_t + x),$$

where $\{W_t\}$ is a standard Brownian motion and $f$ is a bounded function.

Formula (23) makes it possible to solve linear PDE using Monte Carlo methods (the limit law theorem for additive probabilities).

We consider the following example of a nonlinear PDE. Let $u$ be the solution of PDE:

$$\frac{\partial u(t,x)}{\partial t} = \frac{1}{2} \frac{\partial^2 u(t,x)}{\partial x^2} + g\left(u, \frac{\partial u(t,x)}{\partial x}\right)$$

(24)

$$u(0,x) = f(x), \qquad t \geq 0,$$

where $g$ is a function satisfying (H.1)–(H.3) in Section 2.

If there exists a capacity such that the solution of PDE (24) can be represented by a Choquet expectation, then applying the limit law theorem for nonadditive probabilities in Marinacci (1999) would suggest a Monte Carlo-like method could be used to solve nonlinear PDE (24). Unfortunately, our result shows that this is not generally possible.

THEOREM 2. *In the Brownian setting as above, denote by $u_f(t,x)$ the solution of PDE* (24). *If $g(y,z)$ is nonlinear in $z$, then there is no capacity such that the associated Choquet expectation $C$ satisfies $u_f(t,x) = C[f(W_t + x)]$ for all bounded functions $f$ and for all $x$.*

PROOF. Let $\{W_t\}$ be a Brownian motion, by the general Feynman–Kac formula, see, for example, El Karoui, Peng and Quenez (1997) or Ma, Protter and Yong (1994), $u_f(t,x)$ can be represented, by $g$-expectation, that is,

$$u_f(t,x) = \mathcal{E}_g[f(W_t + x)].$$

Applying Theorem 1 and part 1 of Remark 2 completes the proof of this theorem. □

REMARK 3. Theorems 1 and 2 state that if $g$ is nonlinear in $z$, we cannot find a capacity such that the associated Choquet expectation and $g$-expectation $\mathcal{E}_g[f(W_t + x)]$ satisfy

(25)

$$\mathcal{E}_g[f(W_t + x)] = C[f(W_t + x)]$$

for all bounded functions $f$ and for all $x$.

However, if we further restrict $f$ to a set containing only those bounded functions having strictly positive derivatives, we still can find a nonlinear function $g$ and a Choquet expectation such that equation (25) is true. The following is an example.

EXAMPLE 2. Suppose $\mu \neq 0$ is a constant, let $g(z) = \mu|z|$. Obviously $g$ is nonlinear, but we have

$$\mathcal{E}_g[f(W_t + x)] = E_Q[f(W_t + x)]$$

for all bounded functions $f$ with strictly positive derivatives and for all $x$, provided $Q$ is a probability measure defined by

$$(26) \qquad E\left[\frac{dQ}{dP}\Big|\mathcal{F}_T\right] = e^{-(1/2)\mu^2 T + \mu W_T}.$$

Indeed, let $(y_s, z_s)$ be the solution of the BSDE

$$y_s = f(W_t + x) + \int_s^t \mu|z_r|\,dr - \int_s^t z_r\,dW_r, \qquad 0 \le s \le t.$$

By Lemma 8(ii), $z_r > 0$, $r \in [0, t)$ for all bounded functions $f$ with strictly positive derivatives. Thus, the above BSDE is actually a linear BSDE

$$y_s = f(W_t + x) + \int_s^t \mu z_r\,dr - \int_s^t z_r\,dW_r, \qquad 0 \le s \le t.$$

Let $\tilde{W}_r = W_r - \mu r$. Girsanov's lemma then implies that $\{\tilde{W}_r\}$ is a $Q$-Brownian motion under $Q$ denoted in (26). Moreover, the above BSDE can be rewritten as

$$(27) \qquad y_s = f(W_t + x) - \int_s^t z_r\,d\tilde{W}_r, \qquad 0 \le s \le t.$$

Setting conditional expectation $E_Q[\cdot|\mathcal{F}_s]$ on both sides of BSDE (27),

$$y_s = E_Q[f(W_t + x)|\mathcal{F}_s], \qquad 0 \le s \le t.$$

In particular, if we let $s = 0$, by the definition of $\mathcal{E}_g[\cdot]$,

$$\mathcal{E}_g[f(W_t + x)] = y_0.$$

Thus,

$$\mathcal{E}_g[f(W_t + x)] = E_Q[f(W_t + x)]$$

for all bounded functions $f$ with strictly positive derivatives and for all $x$.

Note both that $Q$ does not depend on the choice of $f$ and that mathematical expectation is a Choquet expectation. The proof is complete.

**Acknowledgments.** The work was done while Zengjing Chen was visiting The University of Western Ontario in 2003, whose hospitality he deeply appreciated. All authors thank an anonymous referee for his careful reading of this manuscript—with his help we have been able to significantly improve this work. All remaining errors and omissions remain, of course, our own.

Z. Chen
T. Chen
Department of Mathematics
Shandong University
Jinan 250100
China

M. Davison
Department of Applied Mathematics
University of Western Ontario
London, Ontario
Canada N6A 5B7
e-mail: mdavison@uwo.ca
url: www.apmaths.uwo.ca/˜mdavison